# On Extend the Domain of (Co)convex Polynomial

Malik Saad Al-Muhja[1,2,*]     Amer Himza Almyaly[1]

[1] *Department of Mathematics and Computer Application, College of Sciences, University of Al-Muthanna, Samawa 66001, Iraq*
[2] *Department of Mathematics and statistics, School of Quantitative Sciences, College of Arts and Sciences, Universiti Utara Malaysia, 06010 Sintok, Kedah, Malaysia*

**Abstract.** We will use different way (in this work) from the existing methods in the literature which speaking in the separation of convex sets was carried out by hyperplanes. We are examining the behavior of convex set which is the domain of convex and coconvex polynomial. We simplify this term as (co)convex polynomial herein.
The main goal of the present work is:
What happens if $\mathbb{D}$ is a domain of (co)convex polynomial of $\Delta^{(2)}(Y_s)$, $s \geq 0$ and $x \in \mathbb{D}$? Is $x$ inflection point at $\mathbb{D}$?

**2010 Math. Sub. Classification.** 41A10, 52A30, 52A27

## I. INTRODUCTION

The domain of polynomial is beneficial to provide best approximation to a given function through an approach which stipulates "$f$ is convex function if $\mathrm{epi}(f)$ is a convex subset of $\mathbb{R}^{n+1}$". Subsequently, the expansion of the domain of polynomials is necessary for the assimilation of more than the characteristics of the convex sets such as supporting hyperplane, strongly $\hbar$-hyperplane and separation optimization theorems (see [6, 8, 13]).

The last literature that discussed about best approximation in covex sets involving discrete sets were found in 1979. Ref. [12, 15] introduced the characteristics of approximation to be used in describing second separation theorem. All these characteristics were limited to the approximation theory from an element to a convex set and best approximations by elements of convex sets. But they have been yet to explore that approach in developing approximation which stipulates "complex functions can be approximated by the simpler ones". The new results concept will offer a broader scope for discussion of results in approximation theory, such as an extend the description of separation theorems.

Next, let $X$ be a vector space that has topology $\tau$, then $X$ is locally convex space if any point has a neighborhood base consisting of convex sets (see, [11]). Assumption $\underline{f}$ is continuous convex function from vector space $X$ over field $\mathbb{F}$ onto that field $\mathbb{F}$. Also, the element $g_o \in G$ and

$$\underline{f}(g_o) = \inf\left(\underline{f}(G)\right). \qquad (1)$$

Now, we suppose $\mathbb{X}$ denote the set of all the functions $\underline{f}$ on $X$.
In 1979, Singer [15] proved for any convex subset $G$ of $X$, $(G \neq X)$ satisfying

$$\inf\left(\underline{f}(X)\right) < \inf\left(\underline{f}(G)\right), \qquad (2)$$

---

* Corresponding author's e-mail: dr.al-muhja@hotmail.com & malik@mu.edu.iq

where $X$ is locally convex space. Furthermore, he found

$$\inf\left(\underline{f}(G)\right) = \sup_{\Lambda_1} \inf_{V_1} \underline{f}(y)$$

$$\Lambda_1 = \left\{\underline{f_1} \in \mathbb{X}: \sup\left(\underline{f_1}(G)\right) \leq \underline{f_1}(x), \underline{f_1} \neq 0\right\}$$

$$V_1 = \left\{y \in X: \underline{f_1}(y) = \sup\left(\underline{f_1}(G)\right)\right\}$$

and

$$\inf\left(\underline{f}(G)\right) = \sup_{\Lambda_2} \inf_{V_2} \underline{f}(y)$$

$$\Lambda_2 = \left\{\underline{f_2} \in \mathbb{X}: \sup\left(\underline{f_2}(G)\right) \leq \underline{f_2}(x), \underline{f_2} \neq 0\right\}$$

$$V_2 = \left\{y \in X: \underline{f_2}(y) = \sup\left(\underline{f_2}(G)\right)\right\}$$

where $x \in X$ and

$$\underline{f}(x) < \inf\left(\underline{f}(G)\right). \tag{3}$$

He further proposed some results, if $\underline{f}$ is finite, then (3) is valid.

**Theorem A. [15]** Let $X, \underline{f}, G$ be defined in above, and satisf-ying (2), and $x$ be any element of $X$. If $x$ is satisfying (3), and $g_o \in G$ satisfying (1) iff there exists $\check{f} \in \mathbb{X}, \check{f} \neq 0$, such that

$$\underline{\check{f}}(g_o) = sup\left(\underline{\check{f}}(G)\right) \leq \underline{\check{f}}(x), \text{ and } \underline{f}(g_o) = \inf_{y \in X} \underline{f}(y), \text{ where } \underline{\check{f}}(y) = \underline{\check{f}}(g_o).$$

We will adopt the following concepts in this work.

**Definition 1. [10]** A subset $X$ of $\mathbb{R}^n$ is convex set if $[x, y] \subseteq X$, whenever $x, y \in X$. Equivalently, $X$ is convex if

$$(1 - \lambda)x + \lambda y \in X, \text{ for all } x, y \in X \text{ and } 0 < \lambda < 1.$$

**Definition 2. [5]** The set $\{(x, \mu) \in X \times \mathbb{R}: X \subseteq \mathbb{R}^n, \mu \in \mathbb{R}, \mu \geq f(x)\}$ is called the epigraph of $f$ and denoted by $\text{epi}(f)$. We define the function $f$ on $X$ be a convex function on $X$ if $\text{epi}(f)$ is convex subset of $\mathbb{R}^{n+1}$.

**Remark 3. [14]** A function $f: \mathbb{R}^n \to [-\infty, \infty]$ is convex iff

$$f\big((1 - \lambda)x + \lambda y\big) \leq (1 - \lambda)\alpha + \lambda\beta, 0 < \lambda < 1,$$

whenever $f(x) < \alpha$ and $f(y) < \beta$.

**Theorem 4. [6]** Let $n \in \{2,3\}$ and $f \in \mathbb{C}^{2n}[-1,1]$ be $(2n - 1)$-convex. Then

$$0 \leq \int_{-1}^{1} f(t) - \mathcal{G}_n(f) \leq \mathcal{L}_{n+1}[f] - \int_{-1}^{1} f(t).$$

The following result immediately of strongly $\hbar$-convex (see [3]). In 2016, Lara et al. [8, Corollary 5] proposed a function called $\varepsilon$-strongly $\hbar$-convex such that

$$g(x) - \varepsilon \leq \varphi(x) \leq g(x), x \in D.$$

Let $\pi_n$ be the space of all algebraic polynomials of degree $\leq n - 1$, and $\Delta^{(2)}(Y_s)$ be the collection of all functions $f \in \mathbb{C}[-1,1]$ that change convexity at the points of the set $Y_s$, and are convex in $[y_s, 1]$. The degree of best uniform coconvex polynomial approximation of f is defined by

$$E_n^{(2)}(f, Y_s) = \inf_{p_n \in \pi_n \cap \Delta^{(2)}(Y_s)} \|f - p_n\|$$

where $Y_s = \{y_i\}_{i=1}^{s}$ such that $y_0 = -1 < y_1 < \cdots < y_s < 1 = y_{s+1}$ (see [7]).

If $Y_s = \emptyset$, then $E_n^{(2)}(f, \emptyset) = E_n^{(2)}(f)$ which is usually referred to as the degree of best uniform convex polynomial approximation (see [9]).

**Definition 5.** [4] The weighted Ditzian-Totik moduli of smoothness (DTMS) of $f \in L_p[-1,1]$, when $0 < p \leq \infty$, is defined by

$$\omega_{k,r}^{\phi}(f, t)_p = \sup_{0 < h \leq t} \|\phi(x)^r \Delta_{h\phi}^k(f, x)\|_p$$

where $\phi(x) = \sqrt{1 - x^2}$. If $r = 0$, then

$$\omega_k^{\phi}(f, t)_p = \omega_{k,0}^{\phi}(f, t)_p = \sup_{0 < h \leq t} \|\Delta_{h\phi}^k(f, x)\|_p$$

is the usual DTMS. Also, note that $\omega_{0,r}^{\phi}(f, t)_p = \|\phi^r f\|_p$.

## II. THE MAIN RESULTS

In this section, we will discuss the Domain of Convex Polynomial (DCP). Let $X \subseteq \mathbb{R}$, then

**Definition 1.** A domain $\mathbb{D}$ of convex polynomial $p_n$ of $\Delta^{(2)}$ is a subset of $X$ and $X \subseteq \mathbb{R}$, satisfying the following properties:

1) $\mathbb{D} \in \mathcal{K}^N$, where

$$\mathcal{K}^N = \{\mathbb{D} : \mathbb{D} \text{ is a compact subset of } X\}$$

is the class of all domain of convex polynomial,

2) there is the point $t \in X/\mathbb{D}$, such that

$$|p_n(t)| > \sup\{|p_n(x)| : x \in \mathbb{D}\}, \text{ and}$$

3) there is the function $f$ of $\Delta^{(2)}$, such that

$$\|f - p_n\| \leq \frac{c}{n^2} \omega_{2,2}^{\phi}\left(f'', \frac{1}{2}\right).$$

Let $\mathbb{D}$ and $X$ be as in Definition 1.

**Definition 2.** If the compact set $\mathcal{U}$ is convex, so there is bounded neighborhood set $D = \{\xi \in X : |\xi|^2 < \mathbb{c}\}$ for $\mathbb{c}$ suitably near.

**Theorem 3.** If $\mathbb{D}$ is DCP of $p_n$, and if $x_o \in \mathbb{D}$. Then there is a compact neighborhood $\mathbb{Y}$ of the point $x_o$.

Proof. Suppose that $\mathbb{D}$ is DCP, from Definition 1, then $\mathbb{D}$ is compact subset (CS) of $X$, and $\mathbb{D} \in \mathcal{K}^N$.

Then $\mathbb{D}$ is compact and convex subset of $X$.

From Definition 2, there is bounded neighborhood $\mathbb{Y}$ of the point $x_o$, such that

$$\mathbb{Y} = \{x \in \mathbb{D} : |x|^2 < \mathbb{c}\}, \text{ for } \mathbb{c} \text{ suitably near,}$$

and $\mathbb{Y} \subseteq \mathbb{D}$.

Since $x_o \in \mathbb{D}$, and $|x_o|^2 < \mathbb{c}$, for $\mathbb{c}$ suitably near. Then, $\mathbb{Y} = \mathbb{D}$. Therefore, $\mathbb{Y}$ is compact neighborhood (CNE) of the point $x_o$, and $\mathbb{Y}$ DCP of $p_n$.

**Corollary 4.** If $x_o \in \mathbb{D}$ is DCP of $p_n$. Then $\mathbb{D}$ is CNE of the point $x_o$.

Proof. Clear.

**Theorem 5.** If $p_n : \mathbb{D} \to \mathbb{D}$ is convex polynomial of $\Delta^{(2)}$, and $\mathbb{Y}$ is CS of $\mathbb{D}$. Then $\mathbb{Y}$ is DCP of $p_n$ if and only if $p_n^{-1}(\mathbb{Y})$ is DCP of $p_n$.

Proof. Suppose that $\mathbb{Y}$ is CS of $\mathbb{D}$.

**Case I.** Suppose that $\mathbb{Y}$ is DCP of $p_n$.

Since $p_n: \mathbb{D} \to \mathbb{D}$, and $\mathbb{Y} \subseteq \mathbb{D} \subseteq X$, then $\mathbb{Y}$ is CS of $X$, and $\mathbb{Y} \in \mathcal{K}^N$.

Therefore, $p_n$ is continuous and $p_n^{-1}(\mathbb{Y}) = \mathbb{D}$ is CS of $X$.

Let $t \notin p_n^{-1}(\mathbb{Y})$, then $t \in X/p_n^{-1}(\mathbb{Y})$. From Definition 1, we have
$|p_n(t)| > \sup\{|p_n(x)|: x \in \mathbb{D}\}$, and
the function $f \in \Delta^{(2)}$, such that
$$\|f - p_n\| \leq \frac{c}{n^2} \omega_{2,2}^{\phi}\left(f'', \frac{1}{2}\right).$$
Therefore, $p_n^{-1}(\mathbb{Y})$ is DCP of $p_n$.

**Case II.** Suppose that $p_n^{-1}(\mathbb{Y})$ is DCP of $p_n$.

Since $\mathbb{Y}$ is CS of $\mathbb{D}$.

Let $y \notin \mathbb{D}$, this is, $y \in X/\mathbb{D}$, implies $y \in X/\mathbb{Y}$. From Definition 1, we have
$|p_n(y)| > \sup\{|p_n(x)|: x \in \mathbb{D}\}$,
then,
$|p_n(y)| > \sup\{|p_n(x)|: x \in \mathbb{Y}\}$, and
the function $f \in \Delta^{(2)}$, such that
$$\|f - p_n\| \leq \frac{c_1}{n^2} \omega_{2,2}^{\phi}\left(f'', \frac{1}{2}\right).$$
Therefore, $\mathbb{Y}$ is DCP of $p_n$.

**Theorem 6.** If $\mathbb{D}$ is DCP of $p_n$, and $\mathcal{D} \subseteq \mathbb{D}$ is DCP of $p_n$. For every convex function $f$ of $\Delta^{(2)}$, defined on a neighborhood of $\mathbb{D}$, then the set $\mathcal{D} \cup f^{-1}(0)$ is DCP.

Proof. Suppose that $\mathcal{D} \subseteq \mathbb{D}$ is DCP of $p_n$.

Let $x_o \in \mathcal{D}$, then from Theorem 3, there is CNE $\mathbb{Y}$ of the point $x_o$.

If $f \in \Delta^{(2)}$, such that $f$ define on $\mathbb{Y}$. From Theorem 5, then, $f^{-1}(\mathbb{Y})$ is DCP of $p_n$.

Assume $x_o = 0$, then $\mathcal{D} \cup f^{-1}(0)$ is DCP of $p_n$.

Now, we will define the Domain of Coconvex Polynomial (DCCP).

**Definition 7.** A domain $\mathbb{D}$ of coconvex polynomial $p_n$ of $\Delta^{(2)}(Y_s)$ is a subset of $X$ and $X \subseteq \mathbb{R}$, satisfying the following properties:

1) $\mathbb{D} \in \mathcal{K}^N(Y_s)$, where
$$\mathcal{K}^N(Y_s) = \left\{\begin{matrix}\mathbb{D}: \mathbb{D} \text{ is a compact subset of } X, \\ \text{and } p_n \text{ changes convexity at } \mathbb{D}\end{matrix}\right\}$$
is the class of all domain of coconvex polynomial,

2) $y_i$'s are inflection points, such that
$$|p_n(y_i)| \leq \frac{1}{2}, i = 1, \dots, s, \text{ and}$$

3) there is the function $f$ of $\Delta^{(2)}(Y_s)$, such that
$$\|f - p_n\| \leq \frac{c}{n^2} \omega_{k,2}^{\phi}\left(f'', \frac{1}{n}\right).$$

Let $\mathbb{D}$ and $X$ be as in Definition 7.

**Theorem 8.** If $p_n: \mathbb{D} \to \mathbb{D}$ is coconvex polynomial of $\Delta^{(2)}(Y_s)$ and $\mathbb{D}$ is DCCP of $p_n$. Then $\mathbb{Y}$ is DCCP of $p_n$, if $\mathbb{Y}$ is CNE of the point $x_o$, where $p_n(x_o) = \frac{1}{2}$.

Proof. Suppose that $p_n: \mathbb{D} \to \mathbb{D}$ is coconvex polynomial of $\Delta^{(2)}(Y_s)$, such that $\mathbb{D}$ is CS of $X$, and $p_n$ changes convexity at $\mathbb{D}$.

Put $\mathbb{Y}$ is CNE of $x_o$, implies $x_o \in \mathbb{Y}$.

Since $p_n(x_o) = \frac{1}{2}$, and $\mathbb{D}$ is DCCP of $p_n$. From Definition 7, then:

**Case I.** Either $x_o$ is inflection point at $\mathbb{D}$.

Therefore, $x_o \in \mathbb{D}$, and $\mathbb{Y} \subseteq \mathbb{D}$.

Since, $p_n(x_o) = \frac{1}{2}$. Then, $x_o$ is inflection point at $\mathbb{Y}$.

**Case II.** Or $x_o$ is not inflection point at $\mathbb{D}$.

Now, we must prove that $p_n$ changes convexity at $\mathbb{Y}$. Let $1 \leq s < \infty$, $y_{s-1}, y_s \in \mathbb{D}$, $y_{s-1}, y_s$ be inflection points at $\mathbb{D}$, such that $p_n(y_{s-1}) \leq p_n(x_o) \leq p_n(y_s)$. Since $p_n(x_o) = \frac{1}{2}$, and $y_s$ is inflection points at $\mathbb{D}$, implies $p_n(x_o) = p_n(y_s)$. This is contradiction.

Therefore, $x_o$ is inflection points at $\mathbb{Y}$, and $\mathbb{Y} \subseteq \mathbb{D}$.

Thus, $p_n$ changes convexity at $\mathbb{Y}$.

To prove $\mathbb{Y}$ have all inflection points $\leq \frac{1}{2}$, let $y_j$ be inflection point at $\mathbb{Y}$, such that $j < s$, and $|p_n(y_j)| > \frac{1}{2}$. We get contradiction.

Since $\mathbb{Y} \subseteq \mathbb{D}$, then $f \in \Delta^{(2)}(Y_s)$, such that

$$\|f - p_n\| \leq \frac{c_1}{n^2} \omega_{k,2}^{\phi}\left(f'', \frac{1}{n}\right).$$

Thus, $\mathbb{Y}$ is DCCP of $p_n$.

**Definition 9.** $\gamma$-$H$ is said to be supporting hyperplane to domain of (co)convex polynomial $p_n$ if at least one point $\hat{x}_o$ of $\mathbb{D}$ lies in $\gamma$-$H$, and $p_n(y) \geq \hat{\alpha}$, for all $y \in \mathbb{D} - \{\hat{x}_o\}$, and $\hat{\alpha} \in \mathbb{R}$.

**Definition 10.** If $\mathbb{D}$ is domain of (co)convex polynomial of $p_n$, $x \notin \mathbb{D}$. $\gamma$-$H$ is said to be strictly separates $\mathbb{D}$, if we choose $\mathfrak{b} \in \mathbb{R}$ such that

$$\sup\{p_n(y): y \in \mathbb{D}\} < \mathfrak{b} < p_n(x).$$

**Definition 11.** If $\hbar: [0,1] \to \mathbb{R}$ is given function, $\mathbb{D}_1$ and $\mathbb{D}_2$ are domains of (co)convex polynomials of $p_n$ and $q_n$ respectively. $\gamma$-$H$ and $\gamma_\hbar$-$H$ are said to be strongly hyperplane and strongly $\hbar$-hyperplane respectively, if and only if

$$\inf\{p_n(a): a \in \mathbb{D}_1\} > \sup\{q_n(b): b \in \mathbb{D}_2\}$$

and

$$\inf\{\hbar(t)p_n(a): a \in \mathbb{D}_1\} > \sup\{\hbar(t)q_n(b): b \in \mathbb{D}_2\},$$

where $t \in [0,1]$.

## II. CONCLUSIONS

**Example 1.** Let $n = 3$, $p_3: \mathbb{D} \to (-\infty, \infty)$ be polynomial of degree $\leq n - 1$, such that $\mathbb{D} = [-3,3]$ and $p_3(x) = 0.5x^2 - x$.

1) Suppose that $x = 3$, $y = -3$ and $\lambda = 0.6$ ($0 < \lambda < 1$). Then, $p_3(x = 3) = 1.5$ and $p_3(y = -3) = 7.5$.

Now,

$$p_3((1-\lambda)x + \lambda y) = p_3(-0.6) = -0.78,$$

also,

$$(1-\lambda)p_3(3) + \lambda p_3(-3) = (0.4) \times (1.5) + (0.6) \times (7.5) = 5.1.$$

Therefore,
$$p_3((1 - \lambda) \times (3) + \lambda \times (-3)) \leq (1 - \lambda)p_3(3) + \lambda p_3(-3).$$
Then, $p_3$ is convex polynomial, and $\mathbb{D} \in \mathcal{K}^N$.

**2)** Let $t = 6 \in \mathbb{R}/\mathbb{D}$, then
$$p_3(t = 6) = 12,$$
and
$$\sup\{|p_3(x)|: x \in \mathbb{D}\} = |p_3(x = -3)| = 7.5.$$

**3)** Let $f: \mathbb{D} \to (-\infty, \infty)$ such that
$$f(x) = \begin{cases} \frac{1}{2}x^4 - (x-1)^3 - 2x^2 & ; \text{if } 0 \leq x \leq 3 \\ x & ; \text{if } -3 \leq x < 0 \end{cases},$$
$$f'(x) = \begin{cases} 2x^3 - 3(x-1)^2 - 4x & ; \text{if } 0 \leq x \leq 3 \\ 1 & ; \text{if } -3 \leq x < 0 \end{cases}$$
and
$$f''(x) = \begin{cases} 6x^2 - 6x + 2 & ; \text{if } 0 \leq x \leq 3 \\ 0 & ; \text{if } -3 \leq x < 0 \end{cases}.$$

Let $x_o = 1$, $y_o = 2$ and $\lambda = 0.5$ ($0 < \lambda < 1$). Then, $f(x_o = 1) = -1.5$ and $f(y_o = 2) = -1$.
So,
$$f((1 - \lambda) \times (1) + \lambda(2)) = f(1.5) = -2.093,$$
also,
$$(1 - \lambda)f(1) + \lambda f(2) = -1.25.$$

Therefore,
$$f((1 - \lambda) \times (1) + \lambda(2)) \leq (1 - \lambda)f(1) + \lambda f(2).$$
Hence, $f$ is convex function, and it has $f''$.

Now,
$$\|f(3) - p_3(3)\| = \left\|\left(\frac{1}{2}x^4 - (x-1)^3 - 2x^2\right) - (0.5x^2 - x)\right\| = 13,$$
and
$$\Delta^2_{(0.4)\phi}(f'', x) = \sum_{i=0}^{2} \binom{2}{i}(-1)^{2-i} f''\left(x - \frac{2 \times (0.4)}{2} + i \times (0.4)\right) = 1.92$$

Therefore,
$$\omega^{\phi}_{2,2}\left(f''(x) = 6x^2 - 6x + 2, \frac{1}{2}\right)$$
$$= \sup_{0 < h \leq \frac{1}{2}} \|(1 - x^2) \times \Delta^2_h(f'', x)\|$$
$$= |(-8) \times (1.92)| = 15.36.$$

Thus,
$$\|f - p_3\| \leq \frac{c_1}{9} \omega^{\phi}_{2,2}\left(f''(x) = 6x^2 - 6x + 2, \frac{1}{2}\right),$$
where $c_1 = 7.62$.

**Example 2.** Let $n = 5$, $p_5: \mathbb{D} \to (-\infty, \infty)$ be polynomial of degree $\leq n - 1$, such that $\mathbb{D} = [-3,3]$ and $p_5(x) = (x + 2)(x + 1)(x - 1)(x - 2)$.

**1)** Suppose that $x = 1.5$, $y = 1$ and $\lambda = 0.5$ ($0 < \lambda < 1$). Then, $p_5(x = 1.5) = -2.1875$ and $p_5(y = 1) = 0$.

Now,
$$p_5((1-\lambda)x + \lambda y) = p_5(1.25) = -1.37,$$
also,
$$(1-\lambda)p_5(1.5) + \lambda p_5(1) = (0.5) \times (-2.1875) + (0.5) \times (0) = -1.09.$$
Therefore,
$$p_5((1-\lambda) \times (1.5) + \lambda(1)) \leq (1-\lambda)p_5(1.5) + \lambda p_5(1).$$
Then, $p_5$ is changes convexity at $\mathbb{D} \in \mathcal{K}^N$.

2) Let $Y_s = \{y_i\}_{i=1}^{s=4}$ such that $y_0 = -3 < y_1 = -2 < y_2 = -1 < y_3 = 1 < y_4 = 2 < y_{s+1} = 3$ and are convex in $[y_4, 3]$. Then,
$$|p_5(y_i)| = 0 \leq \frac{1}{2}, i = 1, \dots, 4,$$

3) Let $f: \mathbb{D} \to (-\infty, \infty)$ such that
$$f(x) = \begin{cases} |x^2 - 4| + x & ; \text{if} -3 \leq x \leq 0 \\ |2x - 4| - x & ; \text{if } 0 < x \leq 3 \end{cases},$$

$$f'(x) = \begin{cases} \dfrac{2x^3 - 8x}{|x^2 - 4|} + 1 & ; \text{if} -3 \leq x \leq 0 \\ \dfrac{4x - 8}{|2x - 4|} - 1 & ; \text{if } 0 < x \leq 3 \end{cases}$$

and
$$f''(x) = \begin{cases} \dfrac{(|x^2-4|)^2 \times (6x-8) - (2x^3-8x)^2}{(|x^2-4|)^3} & ; \text{if} -3 \leq x \leq 0 \\ 0 & ; \text{if } 0 < x \leq 3 \end{cases}.$$

Let $x_o = 0, y_o = 0.5$ and $\lambda = 0.5$ ($0 < \lambda < 1$). Then, $f(x_o = 0) = 4$ and $f(y_o = 0.5) = 2.5$. So,
$$f((1-\lambda) \times (0) + \lambda(2.5)) = f(1.25) = 0.25,$$
also,
$$(1-\lambda)f(0) + \lambda f(0.5) = 3.25.$$
Therefore,
$$f((1-\lambda) \times (0) + \lambda(0.5)) \leq (1-\lambda)f(0) + \lambda f(0.5).$$
Hence, $f$ is changes convexity at $\mathbb{D}$, and it has $f''$.

Now,
$$\|f(-3) - p_5(-3)\| = \|(|x^2 - 4| + x) - (x^4 - 5x^2 + 4)\| = 38,$$
and
$$\Delta^4_{(0.1)\phi}(f'', x) = \sum_{i=0}^{4} \binom{4}{i}(-1)^{4-i} f''\left(x - \frac{4 \times (0.1)}{2} + i \times (0.1)\right) = 124.678$$

Therefore,
$$\omega^\phi_{4,2}\left(f''(x) = \frac{(|x^2-4|)^2 \times (6x-8) - (2x^3-8x)^2}{(|x^2-4|)^3}, \frac{1}{5}\right)$$
$$= \sup_{0 < h \leq \frac{1}{2}} \|(1-x^2) \times \Delta^4_{(0.1)}(f'', x)\| = |(-8) \times (124.678)| = 997.424.$$

Thus,
$$\|f - p_5\| \leq \frac{c_2}{25} \omega^\phi_{4,2}\left(f''(x) = \frac{(|x^2-4|)^2 \times (6x-8) - (2x^3-8x)^2}{(|x^2-4|)^3}, \frac{1}{5}\right),$$

where $c_2 = 0.953$.

These results (Definitions 1, 7) are able to answer the question above. Also, it's a possibility of supporting the separation hyperplane theorem later by using DCP (see [1], [2]) like:

If $p_n$ and $q_n$ are two convex polynomials of $\Delta^{(2)}$. If $\mathbb{D}_{p_n}$ is a nonempty compact (and $\mathbb{D}_{q_n}$ is a nonempty closed), such that $\mathbb{D}_{p_n}$ and $\mathbb{D}_{q_n}$ are disjoint. Are $p_n$ and $q_n$ strongly separated by a hyperplane?

## Funding
The research did not receive specific funding yet. The research was performed as part of the employment by the authors.

## Acknowledgements
The authors are indebted to administrative and technical support by University of Al-Muthanna.

## Conflicts of Interest
The authors declare that there is no conflict of interests regarding the publication of this paper.